\documentclass[final,leqno,onetabnum]{siamltex704}

\usepackage{amsmath,amssymb} % for \begin{bmatrix}, \begin{smallmatrix}, \begin{align} etc.
\usepackage{amsfonts} % for \mathbb{} etc. `blackboard' for uppercase only.
\usepackage{exscale}

\usepackage{ifpdf}
\ifpdf
\usepackage[%
  pdftitle={Rayleigh-Ritz majorization error bounds with applications to
FEM},%
  pdfauthor={A. V. Knyazev, M. E. Argentati},%
  pdfsubject={Rayleigh-Ritz method error bounds, FEM analysis},%
  pdfkeywords={Majorization, angles between subspaces,
orthogonal projection, perturbation, a priori error analysis,
Ritz values, Rayleigh-Ritz method,
eigenvalue error bounds,
FEM, finite elements.
},%
  pdfstartview=FitH,%
  bookmarks=true,%
  bookmarksopen=true,%
  breaklinks=true,%
  colorlinks=true,%
  linkcolor=blue,anchorcolor=blue,%
  citecolor=blue,filecolor=blue,%
  menucolor=blue,pagecolor=blue,%
  urlcolor=blue]{hyperref}
\else
\usepackage[%
  breaklinks=true,%
  colorlinks=true,%
  linkcolor=blue,anchorcolor=blue,%
  citecolor=blue,filecolor=blue,%
  menucolor=blue,pagecolor=blue,%
  urlcolor=blue]{hyperref}
\fi

\usepackage{url,doi}

\renewcommand{\ldots}{\ensuremath{\dotsc}}

\newcommand{\comment}[1]{}

\newcommand{\da}{{\downarrow}}
\newcommand{\ua}{{\uparrow}}

\newcommand{\C}{\mathbb{C}} % complex 
\newcommand{\R}{\mathbb{R}} % reals 
\def\Z{{\mathcal Z}}
\def\H{{\mathcal H}}
\newcommand{\X}{\mathcal{X}}
\newcommand{\Y}{\mathcal{Y}}
\newcommand{\XP}{{\mathcal{X}^\perp}}
\newcommand{\YP}{{\mathcal{Y}^\perp}}

\def\lmax{\lambda_{\max}}
\def\lmin{\lambda_{\min}}

\def\size{{{\rm size}}}
\def\span{{\rm span}}
\def\dim{{{\rm dim}}}
\def\max{{{\rm max}}}
\def\min{{{\rm min}}}

\newtheorem{conj}{Conjecture}[section]

\title{Rayleigh-Ritz majorization error bounds with applications to FEM
\thanks{%
Received by the editors May 30, 2008;
accepted for publication (in revised form) September 22, 2009;
this version generated \today.
See also \url{http://arxiv.org/abs/math/0701784}.
}}

\author{%
Andrew V. Knyazev\footnotemark[2]\ \footnotemark[3]\ \footnotemark[4]
\and
Merico E. Argentati\footnotemark[2]\ \footnotemark[3]
}

\begin{document}
\vspace{-1.2in}
%\slugger{simax}{2009}{00}{0}{000-000} %used for `SIAM etc' at the top.
\vspace{.9in}

\setcounter{page}{1}
\maketitle

\renewcommand{\thefootnote}{\fnsymbol{footnote}}
\footnotetext[2]{Department of Mathematical and Statistical Sciences; 
University of Colorado Denver,
P.O. Box 173364, Campus Box 170, Denver, CO 80217-3364, USA.}
\footnotetext[3]{(andrew.knyazev,merico.argentati)[at]ucdenver.edu}
\footnotetext[4]{\url{http://math.ucdenver.edu/~aknyazev/}.
Supported by the NSF-DMS 0612751 and 0728941.}
\renewcommand{\thefootnote}{\arabic{footnote}}

\begin{abstract}
The Rayleigh-Ritz (RR) method finds the stationary values, called Ritz values,
of the Rayleigh quotient on a given trial subspace
as approximations to eigenvalues of a
Hermitian operator $A$. If the trial subspace is $A$-invariant,
the Ritz values are exactly some of the eigenvalues of $A$.
Given two subspaces $\X$ and $\Y$ of the same
finite dimension, such that $\X$ is $A$-invariant,  the
absolute changes in the Ritz values of $A$ with respect to $\X$
compared to the Ritz values with respect to $\Y$ represent the
RR absolute eigenvalue approximation error. 
Our first main result is a sharp  majorization-type RR error bound 
in terms of the principal angles between $\X$ and $\Y$  for an arbitrary
$A$-invariant $\X$, which was a conjecture in 
[SIAM J. Matrix Anal. Appl., 30 (2008), pp. 548-559].
Second, we prove a novel type of RR error bound that deals with the products 
of the errors, rather than the sums. Third, we establish majorization bounds 
for the relative errors. We extend our bounds to the case $\dim\X\leq\dim\Y<\infty$ 
in Hilbert spaces and apply them in the context of the finite element method. 
\end{abstract}

\begin{AM}
15A42, %Inequalities involving eigenvalues and eigenvectors
15A60, %Norms of matrices, numerical range, applications of functional analysis to matrix theory
65F35, %Matrix norms, conditioning, scaling
65N30. %Finite elements, Rayleigh-Ritz and Galerkin methods, finite methods
\end{AM}

%\begin{DOI}
(Place for Digital Object Identifier, to get an idea of the final spacing.)  %10.1137/050630416
%\end{DOI}

\pagestyle{myheadings}
\thispagestyle{plain}
\markboth{ANDREW V. KNYAZEV and MERICO E. ARGENTATI}
{Rayleigh-Ritz majorization error bounds} %50 Characters Limit

\section{Introduction}\label{intro}

For a Hermitian operator $A$ in a finite dimensional inner product space,
the Rayleigh-Ritz method finds the stationary values, called ``Ritz values,''
of the Rayleigh quotient $\lambda(x)=(x,Ax)/(x,x)$  on a given subspace
as approximations to eigenvalues of $A$.
If this ``trial'' subspace is $A$-invariant, i.e.,\ invariant with respect to $A$, 
the Ritz values are exactly some of the eigenvalues of $A$.
Given two finite dimensional subspaces $\X$ and $\Y$ of the same
dimension, such that $\X$ is $A$-invariant,  the
absolute changes in the Ritz values of $A$ with respect to $\X$
compared to the Ritz values with respect to $\Y$ represent the
absolute eigenvalue approximation error.

A priori error bounds for eigenvalues approximated by the Ritz values form one
of the classical subjects in numerical linear algebra and approximation theory. 
Such error bounds are used, e.g.,\ to estimate
convergence rates of iterative methods for matrix eigenvalue problems.
In approximation theory, the  Rayleigh-Ritz method is the most common
technique of approximating eigenvalues and eigenvectors of
operators, e.g.,\ \cite{kvzrs,MR0400004};
and \emph{a priori} error bounds characterize the approximation quality.
Many \emph{a priori} bounds are known; see, e.g.,\ \cite{MR2206452} and references there.

A priori error bounds for eigenvalues in \cite{MR2206452}
are based on the concept of angles between subspaces---one of the major ideas
in multivariate statistics, closely related
to canonical correlations. %, e.g.,\ \cite{hannan}.
This concept also has applications in linear
functional analysis and operator theory. %, e.g.,\ \cite{MR94i:47001,kato}.
The use of angles between subspaces for eigenvalue bounds
is quite natural and may result in elegant and sharp estimates.

Majorization is another classical area of mathematics with numerous applications,
in particular, for estimates involving eigenvalues and singular values.
This paper includes all the necessary material on majorization, but 
should not serve as an introduction to the subject. 
We follow and refer the reader to
\cite{bhatia_book,gk,hj2,mo},
where background and references to original proofs can be found.

In the pioneering results of \cite{dk70}, majorization is applied
to bound eigenvalue errors \emph{a posteriori}
in the framework of angles between subspaces.
A similar approach for \emph{a priori} Rayleigh-Ritz error bounds
is first developed in \cite{akpp06}, e.g.,\ 
a bound with a sharp constant is proved if $\X$ corresponds
to a contiguous set of extreme eigenvalues of $A$. 

Our first major bound, \eqref{eq.hyp} of Theorem \ref{thm_proximity_paper}, 
extends the result of \cite{akpp06} to the general case of an arbitrary
$A$-invariant subspace $\X$, which solves \cite[Conjecture~3.1]{akpp06}.
Moreover, our new proof, with small modifications, 
also covers \eqref{eq.ritz}---the main result of \cite{ka03,ka05}.
Thus, our two bounds \eqref{eq.ritz} and \eqref{eq.hyp} of Theorem \ref{thm_proximity_paper}
supersede all main Rayleigh-Ritz error bounds of \cite{akpp06,ka03,ka05}.  
Our proof is based on a new generalized pinching inequality for singular values
and eigenvalues, Theorem \ref{th:pin}, which is a natural extension of 
the standard pinching  inequality, e.g.,\ \cite[Problem II.5.4]{bhatia_book}.

Next, for the particular case of extreme eigenvalues 
in Theorem \ref{thm_invariant2} we improve bound \eqref{eq.hyp}, 
by replacing  the scalar constant in the bound with a vector of constants 
on the right-hand side. This is a delicate result---if one divides 
both sides of the improved bound, \eqref{eq:SFsinmaj}, by the vector of the constants, 
the statement no longer holds.  

Our second main result, Theorem \ref{thm:mult}, is a majorization Rayleigh-Ritz error bound 
of multiplicative type, which deals with the products of the errors, rather than the sums. 
It allows us to establish majorization bounds for the relative errors, in Theorem~\ref{cor1gr}.
Finally, we extend our bounds to the case $\dim\X\leq\dim\Y<\infty$ 
in infinite dimensional Hilbert spaces, preparing for section \ref{s:FEM}. 

We apply our Rayleigh-Ritz majorization error bounds in the context of the
finite element method (FEM), and briefly show how they improve the constant for a known
FEM eigenvalue error bound from  \cite{MR2206452} in section \ref{s:FEM}.

There are numerous traditional bounds in the form of vector inequalities 
developed over the decades by many mathematicians. 
The question concerning how the traditional and majorization bounds compare naturally arises. 
It appears feasible that properly formulated 
majorization bounds would eventually outperform and thus replace most conventional 
bounds. We already have examples of such a comparison in the present paper, 
but much more work in this direction is needed and will follow. 

\section{Motivation, Conjectures, and Main Results}\label{s:m}
We introduce definitions and after a brief motivation present 
theorems and conjectures on \emph{a priori} majorization
eigenvalue error bounds using principal angles between subspaces.
We describe related results of \cite{akpp06,ka02,ka05,MR2206452}
that precede the developments of the present paper.

\subsection{Basic definitions}
We give only the definition of majorization here and refer the reader to
subsection \ref{sect.majorization} for an overview of some facts
on majorization that we use.
For a real vector $x = \left[ x_{1},\cdots,x_{n} \right]$
let $x^\downarrow$ be the
vector obtained by rearranging the entries of $x$ in an
algebraically decreasing order,
$x_{1}^\downarrow \geq \cdots \geq x_{n}^\downarrow.$
We use the term ``decreasing'' for  ``nonincreasing'',
and ``increasing'' for ``nondecreasing'', for conciseness.
We say that vector $y$ weakly (sub-)majorizes vector 
$x$ and we use the notation $x\prec_wy$ if 
$\sum_{i=1}^k x_i^\downarrow\le\sum_{i=1}^k y_i^\downarrow,\, k=1,\ldots,n$.
If in addition $\sum_{i=1}^n x_i^\downarrow=\sum_{i=1}^n y_i^\downarrow,$  
we write that vector $y$ (strongly) majorizes vector $x$, which is denoted by $x\prec y$.
We overload the notation for scalars, e.g., $1$ may denote the vector of ones.
Nonnegative vectors with different numbers of entries are compared by 
adding zero entries. 

Let $\H$ be a finite dimensional real or complex vector space, equipped
with an inner product $(\cdot,\cdot)$.
%and the corresponding vector and induced operator norms both denoted by $\|\cdot\|$. 
This abstract inner-product space setting is justified by 
generalizing, in subsection \ref{ss:hilb}, some of our results to the 
case of infinite dimensional Hilbert spaces. 
We denote the vector of eigenvalues of a  linear Hermitian operator $A:\H\to\H$
by $\Lambda(A)$, and keep the same notation for Hermitian matrices.
We assume that the vector of eigenvalues $\Lambda$ is
arranged in decreasing order, i.e.,\ $\Lambda=\Lambda^\da$.
Multiple eigenvalues appear in $\Lambda$
repeatedly according to their multiplicities.
We define singular values of a linear operator $B:\H\to\H$ 
as $S(B)=\Lambda\left(\sqrt{B^*B}\right)$, and 
keep the same notation for (rectangular) matrices,  
in which case the operator adjoint $B^*$
is replaced by the complex conjugate matrix transpose $B^H$. 

Let $A:\H \to\H$ be a linear Hermitian operator
and $P_{\X}$ and $P_{\Y}$ be orthoprojectors onto 
subspaces $\X$ and $\Y$ with $\dim\X\leq\dim\Y$.
We first give a brief description of Ritz values.
%Let ${\X}$ be a (so-called ``trial'') subspace of $\H$.
We define the Rayleigh-Ritz operator
$(P_{\X}A)|_{\X}$ on $\X$, where %$P_{\X}$ is the orthogonal projector onto $\X$ and 
$(P_{\X}A)|_{\X}$ denotes the restriction of the
operator $P_{\X}A$ to its invariant subspace $\X$.
%as discussed, e.g., in  \cite[Sect. 11.4, pp. 234--239]{p}.
The eigenvalues $\Lambda\left(\left(P_{\X}A\right)|_{\X}\right)$ are called Ritz values
of the operator $A$ with respect to the subspace $\X$.
In the particular case $\X=\span\{x\}$ for a nonzero vector $x$
we define the Rayleigh quotient
$\lambda(x)=(x,Ax)/(x,x)=\Lambda\left(\left(P_{\X}A\right)|_{\X}\right)$.
If $\X$ is $A$-invariant, the Ritz values  $\Lambda\left(\left(P_{\X}A\right)|_{\X}\right)$  
are some of the eigenvalues of $A.$
For two subspaces $\X$ and $\Y$ of the same
dimension, such that $\X$ is $A$-invariant,  
and $\Y$ approximates $\X$, it is natural to expect that the
Ritz values  $\Lambda\left(\left(P_{\Y}A\right)|_{\Y}\right)$ 
approximate the subset of the eigenvalues of $A$ given by 
$\Lambda\left(\left(P_{\X}A\right)|_{\X}\right)$. 
The absolute changes in the Ritz values of $A$ with respect to $\X$
compared to the Ritz values with respect to $\Y$ thus represent the
absolute eigenvalue approximation error.

%Let $P_{\X}$ and $P_{\Y}$ be orthoprojectors onto 
%subspaces $\X$ and $\Y$ with $\dim\X\leq\dim\Y$.
The vector of cosines squared of principal angles
\emph{from} the subspace $\X$ \emph{to} the subspace $\Y$ is defined by
$\cos^2 \Theta(\X,\Y)= \Lambda^\ua\left(\left(P_{\X}P_{\Y}\right)|_\X\right),$
where the eigenvalues of $(P_{\X}P_{\Y})|_\X$
are rearranged in increasing order, so that
the cosines are increasing, while the angles are defined such that 
$0\leq\Theta=\Theta^\da\leq\pi/2$.
In other words, the cosine squared of angles from $\X$ to $\Y$ are the 
Ritz values of the operator $P_\Y$ on the trial subspace $\X.$
If $\dim\X=\dim\Y$
the definition becomes symmetric with respect to $\X$ and $\Y$,
so it gives the angles \emph{between} subspaces $\X$ and $\Y.$
In the particular case $\X=\span\{x\}$ and $\Y=\span\{y\}$
for unit vectors $x$ and $y$, the vector
$\Theta(\X,\Y)$ has only one component $\theta(x,y)\in[0,\pi/2]$, which is
the acute angle between $x$ and $y$ defined in the standard way, i.e.,\ 
$\cos\theta(x,y)= |(x,y)|$.
If $\H=\C^n$ or $\R^n$, 
let orthonormal columns of matrices $X$ and $Y$ span the
subspaces $\X$ and $\Y$ correspondingly. Then $P_{\X}=XX^H$ and
$P_{\Y}=YY^H$ so
$\left(\cos^2 \Theta(\X,\Y)\right)^\da= \Lambda\left(\left(P_{\X}P_{\Y}\right)|_\X\right)=
\Lambda\left(X^H(YY^H)X\right)=S^2\left(Y^HX\right).$

\subsection{Motivation}\label{sss:m}
Let us first demonstrate that traditional inequalities may not be adequate  
for bounds on Ritz values that involve angles between subspaces. 
Suppose that we bound the vector 
$\left| \Lambda\left((P_{\X}A)|_{\X}\right)-\Lambda\left((P_{\Y}A)|_{\Y}\right)\right|$
of absolute values of matched distances 
between the decreasingly ordered Ritz values
by the vector $\sin\Theta(\X,\Y)$ of the sine 
of principal angles between ${\X}$ and ${\Y}$, 
in order to estimate the influence
of changes in a trial subspace on the Ritz values for the Rayleigh-Ritz method. 

Without majorization, we can compare vectors using component-wise inequalities. 
For example, let $\dim{\X}=\dim{\Y}=2$ and $\sin\Theta(\X,\Y)=[1,0].$
Suppose that we had the inequality 
$\left| \Lambda\left((P_{\X}A)|_{\X}\right)-\Lambda\left((P_{\Y}A)|_{\Y}\right)\right|
\leq C\sin\Theta(\X,\Y)$ with some $A$-dependent constant $C$. 
Since we set the smallest angle between $\X$ and $\Y$ to be zero, 
such an inequality would imply that at least one of the  Ritz values for 
$\X$ is the same as that for $\Y$. This is true only in 
exceptional cases, e.g., if the subspaces $\X$ and $\Y$ intersect 
in an eigenvector of $A$. Thus, changing even a \emph{single} vector  
in the basis of the trial subspace in the Rayleigh-Ritz method 
typically results in changes in \emph{all} Ritz values.
Thus, such an inequality cannot possibly hold. 

Typical known bounds for changes in the Ritz values are inequalities 
only bounding the \emph{largest} change in the  Ritz values by the
\emph{largest} angle between $\X$ and $\Y$. Is it possible to 
take advantage of the knowledge of other angles to get improved bounds for 
the change in the  Ritz values? 
Majorization comes to the rescue as  a natural tool for such bounds. 

\subsection{Sine-based bounds for equidimensional subspaces}
This is the first main result of the paper for trial subspaces of the same dimension.
\begin{theorem} \label{thm_proximity_paper}
Let ${\X}$ and ${\Y}$ be subspaces of
$\H$ with $\dim{\X}=\dim{\Y}$, and let the operator $A$ be Hermitian,
then (see \cite[Remark 4.1]{ka05})
\begin{equation} \label{eq.ritz}
\left| \Lambda\left((P_{\X}A)|_{\X}\right)
-  \Lambda\left((P_{\Y}A)|_{\Y}\right) \right|
\prec_w \left({\lmax}_{(\X+\Y)}-{\lmin}_{(\X+\Y)}\right)~\sin\Theta(\X,\Y),
\end{equation}
where ${\lmin}_{(\X+\Y)}$ and ${\lmax}_{(\X+\Y)}$
are the smallest and largest eigenvalues of the
operator $(P_{\X+\Y}A)|_{\X+\Y}$, respectively.
Also, if one of the subspaces is  $A$-invariant then
\begin{equation} \label{eq.hyp}
\left| \Lambda\left((P_{\X}A)|_{\X}\right)
-  \Lambda\left((P_{\Y}A)|_{\Y}\right) \right|
\prec_w \left({\lmax}_{(\X+\Y)}-{\lmin}_{(\X+\Y)}\right)~\sin^2\Theta(\X,\Y).
\end{equation}
\end{theorem}

If, e.g.,\ $\X$ is $A$-invariant, then $\Lambda\left((P_{\X}A)|_{\X}\right)$
is a subset of eigenvalues of $A$ counting the multiplicities, so
the left-hand side of bound \eqref{eq.hyp}
represents the absolute eigenvalue approximation error in the Rayleigh-Ritz method.

All main proofs for new results here are collected in section \ref{s:p}, e.g.,
we give a unified proof of both bounds of Theorem \ref{thm_proximity_paper}
in section \ref{sec:proximity_paper}. Our proof is shorter and simpler, but 
more sophisticated, compared to that used in \cite{akpp06}, which covers  
a particular case of \eqref{eq.hyp} only. There are two novel key ideas in the proof.  
First, we concatenate
the absolute values in the left hand side of \eqref{eq.ritz} 
(or \eqref{eq.hyp}) with the same values but with the negative sign. 
Second, our new generalized pinching inequality \eqref{eq:genpinching.strong}
of Theorem \ref{thm:genpinching} accurately bounds the concatenated vector. 

Bound \eqref{eq.ritz} is proved in  \cite{ka05} and bound \eqref{eq.hyp}
is conjectured in \cite[Conjecture~3.1]{akpp06}, but both with a larger constant,
which is the spread  $\lmax-\lmin$ of the spectrum of $A$ 
where $\lmin$ and $\lmax$ are the smallest and largest eigenvalues of $A$, respectively.
Bounds \eqref{eq.ritz} and \eqref{eq.hyp} use 
the smaller spread, ${\lmax}_{(\X+\Y)}-{\lmin}_{(\X+\Y)}\leq\lmax-\lmin$, of the spectrum
of the operator $(P_{\X+\Y}A)|_{\X+\Y}$. However, \cite[Remark 4.1]{ka05}
states that these two statements are in fact equivalent.
The argument of \cite[Remark 4.1]{ka05} is essential
and used several times below. For completeness, let us reproduce it here.

In this paper, we always assume that both
subspaces $\X$ and $\Y$ are finite dimensional, cf.\ \cite{kja07}.
Let us consider the
finite dimensional subspace $\X+\Y$ and the operator $(P_{\X+\Y}A)|_{\X+\Y}$
as replacements for the original space $\H$ and the original operator $A$.
Whether we define the Rayleigh-Ritz operator, e.g.,\
$(P_{\X}A)|_{\X}$ on $\X$ starting with the original
space $\H$ and the operator $A$, or with the reduced
space $\X+\Y$ and the operator $(P_{\X+\Y}A)|_{\X+\Y}$,
the outcome is evidently the same:
$(P_{\X}A)|_{\X}=(P_{\X}(P_{\X+\Y}A)|_{\X+\Y})|_{\X}.$

Moreover,  if $\X$ is $A$-invariant, it is also
$(P_{\X+\Y}A)|_{\X+\Y}$-invariant,
corresponding to the same set of eigenvalues.
If this set of eigenvalues is
the contiguous set of the largest eigenvalues of $A$,
it is also the contiguous set of the largest eigenvalues of $(P_{\X+\Y}A)|_{\X+\Y}$.
Thus, without loss of generality, we can substitute
the space $\X+\Y$ and the operator $(P_{\X+\Y}A)|_{\X+\Y}$
for the space $\H$ and the operator $A$. 
This simple substitution improves the constants and, most importantly, allows us
to handle easily the case of infinite dimensional $\H$
as we explain later in Lemma \ref{l:csle}.

\subsection{Improved sine-based bounds for equidimensional subspaces}\label{ss:imp}
In the right-hand sides of both bounds in Theorem \ref{thm_proximity_paper}  the scalar 
${\lmax}_{(\X+\Y)}-{\lmin}_{(\X+\Y)}$ appears. We want to improve 
Theorem \ref{thm_proximity_paper} by replacing the scalar factor
with the following decreasing \emph{vector} of different scalar factors: 
\[{\rm Spr}_{(\X+\Y)}=
\left[{\lambda_i}_{(\X+\Y)}-{\lambda_{-i}}_{(\X+\Y)},\,i=1,\ldots,\dim\X\right],\]
%is the nonnegative decreasing vector with the first $\dim(\X+\Y)-\dim\Y$ components
%given by ${\lambda_i}_{(\X+\Y)}-{\lambda_{-i}}_{(\X+\Y)},$ 
where ${\lambda_1}_{(\X+\Y)}\geq\cdots\geq{\lambda_{\dim\X}}_{(\X+\Y)}$
and ${\lambda_{-1}}_{(\X+\Y)}\leq\cdots\leq{\lambda_{-\dim\X}}_{(\X+\Y)}$
are the $\dim\X$ largest and smallest, respectively, eigenvalues of the
operator $(P_{\X+\Y}A)|_{\X+\Y}$. Since there are  
$\dim(\X+\Y)-\dim\Y$ nonzero components in the vector 
$\Theta(\X,\Y)$, only the first $\dim(\X+\Y)-\dim\Y$
components of ${\rm Spr}_{(\X+\Y)}$, which are all nonnegative, 
will actually be used in upcoming bounds  \eqref{eq.hyp1} and \eqref{eq.hyp2}, 
where the component-wise products of vectors ${\rm Spr}_{(\X+\Y)}$ 
and $\sin\Theta(\X,\Y)$ appear in the right-hand sides.

Let us consider an extreme case example with $\sin\Theta(\X,\Y)=[1,1]$, i.e.,\ 
uncorrelated $\X$ and $\Y$. The largest variation in the individual Ritz values  
is in this case clearly bounded by the scalar value ${\lmax}_{(\X+\Y)}-{\lmin}_{(\X+\Y)}$
of the spread of the spectrum of  $(P_{\X+\Y}A)|_{\X+\Y}$, 
which is already used in Theorem \ref{thm_proximity_paper} and which is the first 
component in the vector ${\rm Spr}_{(\X+\Y)}$. Now let us consider the sum 
of both components of the vector 
$\left|\Lambda\left((P_{\X}A)|_{\X}\right)-\Lambda\left((P_{\Y}A)|_{\Y}\right)\right|$.
This sum takes the largest value if $\X$ is the span of two eigenvectors of 
$A$ corresponding to its largest eigenvalues, 
while $\Y$ is the span of two eigenvectors of 
$A$ corresponding to its smallest eigenvalues. 
But this largest value in this example is 
exactly the sum of both components of the vector  ${\rm Spr}_{(\X+\Y)}$. 
This example suggests that the  vector  spread ${\rm Spr}_{(\X+\Y)}$ might be 
the appropriate vector of constants to replace the scalar spread  
${\lmax}_{(\X+\Y)}-{\lmin}_{(\X+\Y)}$ in Theorem \ref{thm_proximity_paper}. 
Our numerical tests motivate the following conjecture. 
\begin{conj} \label{hyp2}
Let ${\X}$ and ${\Y}$ be subspaces of
$\H$ with $\dim{\X}=\dim{\Y}$ and operator $A$ be Hermitian.
Then
\begin{equation} \label{eq.hyp1}
\left| \Lambda\left((P_{\X}A)|_{\X}\right)
-  \Lambda\left((P_{\Y}A)|_{\Y}\right) \right|
\prec_w {\rm Spr}_{(\X+\Y)} \sin\Theta(\X,\Y).
\end{equation}
If in addition one of the subspaces is $A$-invariant then
\begin{equation} \label{eq.hyp2}
\left| \Lambda\left((P_{\X}A)|_{\X}\right)
-  \Lambda\left((P_{\Y}A)|_{\Y}\right) \right|
\prec_w {\rm Spr}_{(\X+\Y)} \sin^2\Theta(\X,\Y).
\end{equation}
\end{conj}
Since
$\max\left\{{\rm Spr}_{(\X+\Y)}\right\}={\lmax}_{(\X+\Y)}-{\lmin}_{(\X+\Y)}$,
bounds \eqref{eq.ritz} and \eqref{eq.hyp} would 
follow from \eqref{eq.hyp1} and \eqref{eq.hyp2}, correspondingly.

Here, we are able to prove only \eqref{eq.hyp2}
and under an additional assumption.
\begin{theorem} \label{thm_invariant2}
If $\X$ is $A$-invariant and 
corresponds to the contiguous set of the largest eigenvalues of $A$, then 
bound \eqref{eq.hyp2} holds. Consequently, we obtain  
\begin{equation}
0 \leq
\Lambda\left((P_{\X}A)|_{\X}\right)-\Lambda\left((P_{\Y}A)|_{\Y}\right)
\prec_w
\left(\Lambda\left((P_{\X}A)|_{\X}\right) - {\lmin}_{(\X+\Y)}\right)
\sin^2 \Theta(\X,\Y),
\label{eq:SFsinmaj}
\end{equation}
where we take into account that in this case
${\rm Spr}_{(\X+\Y)}\leq\Lambda\left((P_{\X}A)|_{\X}\right)-{\lmin}_{(\X+\Y)}.$
\end{theorem}

Bounds  \eqref{eq.hyp1}, \eqref{eq.hyp2}, and \eqref{eq:SFsinmaj} are delicate. 
An attempt to improve them by dividing both sides by the vector of the constants breaks them all,
as can be checked by running the following MATLAB code (see \cite{ka02} for SUBSPACEA's description):
\begin{verbatim}
A=diag([2 1 0 0]); X=[1 0;0 1;0 0;0 0]; Y=orth([1 0;1 2;2 -2;0 1]);
Spr=[2 1]'; SinTheta=flipud(sort(sin(subspacea(X,Y))))
LeftHandSide=flipud(sort(abs(sort(eig(X'*A*X))-sort(eig(Y'*A*Y)))))
sum(LeftHandSide)<=sum(Spr.*SinTheta.*SinTheta)   %(2.4) and (2.5) hold 
sum(LeftHandSide./Spr)<=sum(SinTheta)             %divided (2.3) breaks 
\end{verbatim}
In this example, the  vector of the constants in \eqref{eq:SFsinmaj} 
is equal to ${\rm Spr}_{(\X+\Y)}$, 
i.e.,\ bounds  \eqref{eq.hyp2} and  \eqref{eq:SFsinmaj} are the same and must hold, since 
the assumptions of Theorem \ref{thm_invariant2} are satisfied, 
which is confirmed by the code. 
The last line shows that even \eqref{eq.hyp1}, the 
weakest of the three bounds, does not hold if divided by ${\rm Spr}_{(\X+\Y)}$. 

Our bound \eqref{eq:SFsinmaj}  is competitive compared to the following known inequality,
\begin{eqnarray}
\nonumber 0&\leq&
\Lambda\left((P_{\X}A)|_{\X}\right)-\Lambda\left((P_{\Y}A)|_{\Y}\right)\\
&\leq& \left(\Lambda\left((P_{\X}A)|_{\X}\right) - {\lmin}_{(\X+\Y)}\right)
\max\left\{\sin^2 \Theta(\X,\Y)\right\}), \label{eq:SF}
\end{eqnarray}
which is proved in \cite{k86} and presented above in a slightly modified formulation to make it
consistent with \eqref{eq:SFsinmaj}.
There is no majorization in \eqref{eq:SF}---each vector component is bounded separately, thus one 
can divide  \eqref{eq:SF} by the vector of the constants in contrast to \eqref{eq:SFsinmaj}. 
But \eqref{eq:SF} only uses the largest angle, so it would be inferior to \eqref{eq:SFsinmaj} if 
other angles are much smaller compared to the largest angle---a common situation in applications;
e.g.,\ see section \ref{s:FEM} on FEM. 

\subsection{Multiplicative and tangent-based bounds}\label{ssec:mult}
The main goal of this subsection is to formulate \emph{multiplicative} analogs for the 
majorization-type bounds of the previous subsection, 
based on products rather than sums. 

The definition of majorization for vectors is based on the \emph{sums} of 
vector components. It can also deal with the \emph{products} of 
nonnegative vector components using the following conventions. 
If for nonnegative vectors  $x = \left[ x_{1},\cdots,x_{n} \right]$
and  $y = \left[ y_{1},\cdots,y_{n} \right]$ it holds that 
$\prod_{i=1}^k x_i^\downarrow\le\prod_{i=1}^k y_i^\downarrow,\, k=1,\ldots,n$,
we write $\log x\prec_w\log y$.
If in addition $\prod_{i=1}^n x_i^\downarrow=\prod_{i=1}^n y_i^\downarrow$  
we write  $\log x\prec\log y$. For strictly positive vectors these 
conventions follow directly from the definition of (weak) majorization. 

The example of subsection \ref{sss:m} implies that it is impossible to 
bound products of changes of different Ritz values
by the products of the sine of changes in the principal angles. 
Our novel multiplicative bound below uses $1+\tan^{2}$ rather than $\sin^2$ 
to bound the relative eigenvalue error in the form of the products. 
\begin{theorem} \label{thm:mult}
Under the assumptions of Theorem \ref{thm_invariant2}
let $\Theta(\X,\Y)<\pi/2$ 
and $\Lambda\left((P_{\X}A)|_{\X}\right)>{\lmin}_{(\X+\Y)}$. 
Then 
$\Lambda\left((P_{\Y}A)|_{\Y}\right)>{\lmin}_{(\X+\Y)}$ and we have 
\begin{equation*}
0\leq
\log\frac{\Lambda\left((P_{\X}A)|_{\X}\right)-{\lmin}_{(\X+\Y)}}
{\Lambda\left((P_{\Y}A)|_{\Y}\right)-{\lmin}_{(\X+\Y)}}
\prec_w\log\left(1+\tan^{2}\Theta(\X,\Y)\right),
\end{equation*}
which leads to 
\begin{equation*}
0 \leq
\frac{\Lambda\left((P_{\X}A)|_{\X}\right)-\Lambda\left((P_{\Y}A)|_{\Y}\right)}
{\Lambda\left((P_{\Y}A)|_{\Y}\right) - {\lmin}_{(\X+\Y)}}
\prec_w\tan^2 \Theta(\X,\Y).
\end{equation*}
\end{theorem}

We note that either majorization result of Theorem \ref{thm:mult} implies the bound 
\[
0 \leq
\frac{\Lambda\left((P_{\X}A)|_{\X}\right)-\Lambda\left((P_{\Y}A)|_{\Y}\right)}
{\Lambda\left((P_{\Y}A)|_{\Y}\right) - {\lmin}_{(\X+\Y)}}
\leq\max\left\{\tan^2 \Theta(\X,\Y)\right\},
\]
which is an equivalent form of the already known sine-based inequality \eqref{eq:SF}. 

\subsection{Bounds for non-equidimensional subspaces}
In all statements we have made so far we have assumed that $\dim{\X}=\dim{\Y},$
but applications require the more general assumption $\dim{\X}\leq\dim{\Y}$,
where we interpret the principal angles $\Theta(\X,\Y)$ as the angles 
\emph{from} $\X$ \emph{to} $\Y$.
In this paper, we briefly consider one such well known application of 
the Rayleigh-Ritz method: the finite element method for partial differential equations, 
in section \ref{s:FEM}.
%, and subspace iterations, in section \ref{s:si}.

Since we compare  $\dim{\X}\leq\dim{\Y}$ Ritz values for the trial subspace $\X$ 
against $\dim{\Y}$ Ritz values for the trial subspace $\Y$, 
we can either specifically choose some appropriate $\dim{\X}$ Ritz values out of 
 $\dim{\Y}$ Ritz values for $\Y$, or simply state that \emph{there exist}  
$\dim{\X}$ Ritz values for the trial subspace $\Y$ such that our bounds hold. 

Bounds  \eqref{eq.ritz} and \eqref{eq.hyp1} do not hold if 
$\dim{\X}<\dim{\Y}$ even using the latter, weaker, statement. 
Indeed, e.g.,\ if $\dim{\X}=1$ and $\Theta(\X,\Y)=0$ then 
either bound  \eqref{eq.ritz} or \eqref{eq.hyp1} would imply that 
$\Lambda\left((P_{\X}A)|_{\X}\right)$---in this case
a single number---is one of the Ritz values for the trial subspace $\Y$, 
which is not true since $\X$ is arbitrary in $\Y$. 

Known results, e.g.,\ \cite{k97,MR2206452},
guarantee the existence of $\dim{\X}$ Ritz values for the trial subspace $\Y$ 
that are good approximations for $\dim{\X}<\dim{\Y}$ eigenvalues 
for an arbitrary $A$-invariant subspace $\X$ if $\Theta(\X,\Y)$ is small. 
However, our numerical tests show that \eqref{eq.hyp} and \eqref{eq.hyp2} 
still fail in this case; cf. \cite[Lemma 2.6]{MR2206452}.    
An approach of \cite[Theorem 2.7]{MR2206452} may help to overcome 
the obstacle, but it is outside of the scope of this paper.
Here we consider only the particular case where the $A$-invariant subspace $\X$ 
corresponds to the contiguous set of the largest eigenvalues
$\Lambda\left((P_{\X}A)|_{\X}\right)$ of $A$.
%We approximate $\Lambda\left((P_{\X}A)|_{\X}\right)$ using
%the $\dim{\X}\leq\dim{\Y}$ largest Ritz
%values---the eigenvalues of $(P_{\Y}A)|_{\Y}$
%---denoted by $\Lambda_{\dim{\X}}\left((P_{\Y}A)|_{\Y}\right)$.
\begin{theorem} \label{cor1g}
Let $\dim{\X}\leq\dim{\Y}$, the operator $A$ be Hermitian,
the $A$-invariant subspace $\X$
correspond to the contiguous set of the largest eigenvalues of $A$,
and $\Lambda_{\dim{\X}}\left((P_{\Y}A)|_{\Y}\right)$ denote
the $\dim\X$ %(counting the multiplicities) 
largest eigenvalues of $(P_{\Y}A)|_{\Y}.$ Then
\begin{eqnarray}\label{eq:1}
0&\leq&
\Lambda\left((P_{\X}A)|_{\X}\right)-\Lambda_{\dim{\X}}\left((P_{\Y}A)|_{\Y}\right)\\ \nonumber
&\prec_w&
\left(\Lambda\left((P_{\X}A)|_{\X}\right) - {\lmin}_{(\X+\Y)}\right)
\sin^2 \Theta(\X,\Y);
%\label{eq:SFsinmajg}
\end{eqnarray}
if $\Theta(\X,\Y)<\pi/2$ and $\Lambda\left((P_{\X}A)|_{\X}\right)>{\lmin}_{(\X+\Y)}$, then
$\Lambda_{\dim{\X}}\left((P_{\Y}A)|_{\Y}\right)>{\lmin}_{(\X+\Y)},$  
\begin{equation}\label{eq:2}
0\leq
\log\frac{\Lambda\left((P_{\X}A)|_{\X}\right)-{\lmin}_{(\X+\Y)}}
{\Lambda_{\dim{\X}}\left((P_{\Y}A)|_{\Y}\right)-{\lmin}_{(\X+\Y)}}
\prec_w\log\left(1+\tan^{2}\Theta(\X,\Y)\right),
\end{equation}
and
\begin{equation}\label{eq:3}
0 \leq
\frac{\Lambda\left((P_{\X}A)|_{\X}\right)-\Lambda_{\dim{\X}}\left((P_{\Y}A)|_{\Y}\right)}
{\Lambda_{\dim{\X}}\left((P_{\Y}A)|_{\Y}\right) - {\lmin}_{(\X+\Y)}}
\prec_w\tan^2\Theta(\X,\Y).
\end{equation}
\end{theorem}
\begin{proof}
We use a technique presented in \cite{k86}
to extend Rayleigh-Ritz error bounds for the particular case
$\dim{\X}=\dim{\Y}$ to the general case $\dim{\X}\leq\dim{\Y}$.

Let $\dim{\X}\leq\dim{\Y}$ and $\Theta(\X,\Y)<\pi/2.$
We define a new subspace $\Z$ to be the orthogonal projection of $\X$ onto $\Y$,
i.e.,\ $\Z=P_\Y\X$. Assuming $\Theta(\X,\Y)<\pi/2$ gives $\dim\X=\dim\Z$
and $\Theta(\X,\Y)=\Theta(\X,\Z).$
Since $\Z\subseteq\Y,$
the Courant-Fisher min-max principle evidently implies that
$\Lambda\left((P_{\Z}A)|_{\Z}\right)\leq
\Lambda_{\dim{\X}}\left((P_{\Y}A)|_{\Y}\right)$
for the largest $\dim\X=\dim\Z$ eigenvalues and that ${\lmin}_{(\X+\Z)}\geq{\lmin}_{(\X+\Y)},$
so \eqref{eq:SFsinmaj} leads to \eqref{eq:1}.
Since bound \eqref{eq:1} depends continuously on  $\Theta(\X,\Y)$,  
the assumption  $\Theta(\X,\Y)<\pi/2$ can be removed by the continuity argument. 

Now we apply \eqref{eq:SF} to the pair of subspaces $\X$ and $\Z$
instead of $\X$ and $\Y$, i.e.,\
\begin{eqnarray*}
\nonumber 0&\leq&
\Lambda\left((P_{\X}A)|_{\X}\right)-\Lambda_{\dim{\X}}\left((P_{\Y}A)|_{\Y}\right)\leq
\Lambda\left((P_{\X}A)|_{\X}\right)-\Lambda\left((P_{\Z}A)|_{\Z}\right)\\
&\leq& \left(\Lambda\left((P_{\X}A)|_{\X}\right) - {\lmin}_{(\X+\Z)}\right)
\max\left\{\sin^2 \Theta(\X,\Z)\right\})\\
&\leq& \left(\Lambda\left((P_{\X}A)|_{\X}\right) - {\lmin}_{(\X+\Y)}\right)
\max\left\{\sin^2 \Theta(\X,\Y)\right\}).
\end{eqnarray*}
This gives the following known inequality, e.g.,\ \cite{k86},
\begin{eqnarray}\label{eq:SFgn}
0&\leq&
\Lambda\left((P_{\X}A)|_{\X}\right)-\Lambda_{\dim{\X}}\left((P_{\Y}A)|_{\Y}\right)\\ \nonumber
&\leq&\left(\Lambda\left((P_{\X}A)|_{\X}\right) - {\lmin}_{(\X+\Y)}\right)
\max\left\{\sin^2\Theta(\X,\Y)\right\}.\end{eqnarray}
If $\Theta(\X,\Y)<\pi/2$ and $\Lambda\left((P_{\X}A)|_{\X}\right)>{\lmin}_{(\X+\Y)}$, then
$\Lambda_{\dim{\X}}\left((P_{\Y}A)|_{\Y}\right)>{\lmin}_{(\X+\Y)}.$  
Theorem \ref{thm:mult}  
immediately leads to \eqref{eq:2} and \eqref{eq:3} by monotonicity arguments.  
\end{proof}

\subsection{Relative eigenvalue error bounds}
Our previous results bound the absolute value of the eigenvalue error. 
They are all invariant with respect to shifting the  
operator $A$ into $A+\alpha I$ for any real shift $\alpha$.
For eigenvalues that are small in absolute value, it is also important to 
bound the \emph{relative} error. Here we show how new relative 
bounds can be easily obtained from our Theorem \ref{cor1g}. 
For relative bounds the shift-invariance will of course 
be lost, and it is natural to assume that $A>0$. 

Let us first explain how relative bounds are obtained, e.g.,\ from \eqref{eq:SFgn}. 
Since  $A>0$ we can bound ${\lmin}_{(\X+\Y)}\geq0$  
and divide both sides of the inequality by the vector 
$\Lambda\left((P_{\X}A)|_{\X}\right)>0$, which gives 
\begin{equation}\label{e:frac}
0\leq1-\frac
{\Lambda_{\dim{\X}}\left((P_{\Y}A)|_{\Y}\right)}{\Lambda\left((P_{\X}A)|_{\X}\right)}
\leq\max\left\{\sin^2\Theta(\X,\Y)\right\}.
\end{equation}
This is already a relative bound, but only for the largest eigenvalues, which is not so useful.
We can turn the largest eigenvalues into the 
smallest ones by substituting $A^{-1}$ for $A$ as $A>0$, but this substitution 
alone does not reproduce the inverse of the  Rayleigh quotient since in general 
$(x,Ax)(x,A^{-1}x)\neq1.$ There is a simple fix, though.
Introducing the notation 
$(x,y)_A=(x,Ay)$ for the $A$-based scalar product, we have the following 
trivial but crucial identity for the Rayleigh quotient, 
\[
\frac{(x,Ax)}{(x,x)}=\frac{(x,x)_A}{(x,A^{-1}x)_A}=
\left(\frac{(x,A^{-1}x)_A}{(x,x)_A}\right)^{-1}.
\] 
It implies that the Rayleigh-Ritz method on a trial subspace $\X$ 
applied to the operator $A$ in the original scalar product $(\cdot,\cdot)$
or to the operator $A^{-1}$ in the $A$-based scalar product $(\cdot,\cdot)_A$
gives the same Ritz vectors, and the Ritz values are reciprocals of each other. 
The use of the $A$-based scalar product changes the way we measure the angles, 
see~\cite{ka02}. Simultaneous substitutions $A^{-1}$ for $A$ and 
$(\cdot,\cdot)_A$ for $(\cdot,\cdot)$
in \eqref{eq:SFgn} and Theorem \ref{cor1g} give the following new relative bounds. 
\begin{theorem} \label{cor1gr}
Let $\dim{\X}\leq\dim{\Y}$ and $\Theta(\X,\Y)<\pi/2,$ the operator $A$ be Hermitian
and positive definite, $A>0$,
the $A$-invariant subspace $X$
correspond to the contiguous set of the smallest eigenvalues of $A$,
$\Lambda_{\dim{\X}}\left((P_{\Y}A)|_{\Y}\right)$ denote
the $\dim\X$ (counting the multiplicities) smallest
eigenvalues of $(P_{\Y}A)|_{\Y},$ and $\Theta_A(\X,\Y)$ denote the 
vector of angles from $\X$ to $\Y$ defined in the  $A$-based scalar product 
 $(\cdot,\cdot)_A$. Then
\begin{eqnarray}\label{eq:SFg}
\quad 0&\leq&1-\frac
{\Lambda\left((P_{\X}A)|_{\X}\right)}{\Lambda_{\dim{\X}}\left((P_{\Y}A)|_{\Y}\right)}
\leq\max\left\{\sin^2\Theta_A(\X,\Y)\right\},\\\label{e:cos}
\quad 0&\leq&
\log\frac{\Lambda_{\dim{\X}}\left((P_{\Y}A)|_{\Y}\right)}
{\Lambda\left((P_{\X}A)|_{\X}\right)}
\prec_w\log\left(1+\tan^{2}\Theta_A(\X,\Y)\right),\\\label{eq:SFsinmajginf}
\quad 0&\leq&
\frac{\Lambda_{\dim{\X}}\left((P_{\Y}A)|_{\Y}\right)}
{\Lambda\left((P_{\X}A)|_{\X}\right)}
-1
\prec_w\tan^2\Theta_A(\X,\Y).
\end{eqnarray}
\end{theorem}

Let us highlight that bound \eqref{e:cos} is not only relative but 
also multiplicative. 

We finally note that the first statement, with the sine, in Theorem \ref{cor1g} 
cannot be transformed into a relative bound in the same way. A seemingly natural  
extension 
$1-{\Lambda\left((P_{\X}A)|_{\X}\right)}/{\Lambda_{\dim{\X}}\left((P_{\Y}A)|_{\Y}\right)}
\prec_w\sin^2\Theta_A(\X,\Y)$
of \eqref{eq:SFg} is in fact wrong; see 
\begin{verbatim}
A=diag([1 2 3 100]);X=[1 0;0 1;0 0;0 0];Y=orth([-6 -1;-7 1;2 6;1 -7]); 
SinThetaA=flipud(sort(sin(subspacea(X,Y,A))));
LeftHandSide=[1 1]'-[2 1]'./flipud(sort(eig(Y'*A*Y)));
sum(LeftHandSide)<=sum(SinThetaA.*SinThetaA)             %fails 
\end{verbatim}

\subsection{Generalizations for Hilbert spaces}\label{ss:hilb}
Here we extend some of the previous results to
infinite dimensional spaces, using again \cite[Remark 4.1]{ka03}. 
Let $\H$ be an \emph{infinite dimensional} Hilbert space and
$A:\H \to\H$ be a linear bounded Hermitian operator.
Let $P_{\X}$ and $P_{\Y}$ be orthogonal projectors onto the nontrivial
\emph{finite dimensional} subspaces $\X$ and $\Y$ with $\dim\X\leq\dim\Y<\infty$.
The vector of cosines squared of $\dim\X$
principal angles from  $\X$ to $\Y$ is defined by
$\cos^2 \Theta(\X,\Y)= \Lambda((P_{\X}P_{\Y})|_\X).$
If $\X$ is $A$-invariant,
the Ritz values  $\Lambda\left(\left(P_{\X}A\right)|_{\X}\right)$
are some of the eigenvalues of $A$, since we assume that $\X$ is finite dimensional.
Throughout the section, we use the vectors of
eigenvalues $\Lambda$ enumerated in decreasing order
only for finite dimensional operators, so the vectors
have a finite number of components as before.

Both subspaces $\X$ and $\Y$ are finite dimensional. Let us consider the
finite dimensional subspace $\X+\Y$ and the operator $(P_{\X+\Y}A)|_{\X+\Y}$
as replacements to the original space $\H$ and the operator $A$.
The Rayleigh-Ritz operator $(P_{\X}A)|_{\X}$ on $\X$ using the original
space $\H$ and the operator $A$ is the same as using the reduced
space $\X+\Y$ and the operator $(P_{\X+\Y}A)|_{\X+\Y}$,
as we have already discussed.

If $\X$ is $A$-invariant, it is also $(P_{\X+\Y}A)|_{\X+\Y}$-invariant,
corresponding to the same set of eigenvalues.
If this set of eigenvalues is
the contiguous set of the largest eigenvalues of $A$, which
forms the top of the spectrum of $A$, then
it is also the contiguous set of the largest eigenvalues of $(P_{\X+\Y}A)|_{\X+\Y}$.
The latter may not be so evident in the infinite dimensional setting,
so let us give and prove here the formal statement.
\begin{lemma}\label{l:csle}
For a linear bounded Hermitian operator $A$ on an
infinite dimensional Hilbert space $\H$, let $\X$ be
a nontrivial finite dimensional $A$-invariant subspace of $\H$
that corresponds to the top part of the spectrum of $A$,
i.e.,\ the smallest point of the spectrum of
$(P_{\X}A)|_{\X}$ is an upper bound for the
largest point of the spectrum of $(P_{\XP}A)|_{\XP}$.
Then for any nontrivial finite dimensional subspace $\Y$ of $\H$
the Hermitian operator $(P_{\X+\Y}A)|_{\X+\Y}$ is invariant on $\X$,
and the spectrum of the restriction of $(P_{\X+\Y}A)|_{\X+\Y}$
to $\X$ comprises the $\dim\X$ largest eigenvalues of $(P_{\X+\Y}A)|_{\X+\Y}$.
\end{lemma}
\begin{proof}
The spectrum of a bounded Hermitian operator is a closed bounded set on the real line.
Since $\X$ is finite dimensional, the spectrum
$\Lambda\left((P_{\X}A)|_{\X}\right)$ consists of $\dim\X$ eigenvalues, counting the
multiplicities. Since $\X$ is $A$-invariant, the spectrum
$\Lambda\left((P_{\X}A)|_{\X}\right)$ is a subset of the spectrum of $A$, which
by the lemma assumption forms the top part of the spectrum of $A$.
The subspace $\X$ is $A$-invariant by assumption
and is evidently $P_{\X+\Y}$-invariant, so it is also
$(P_{\X+\Y}A)|_{\X+\Y}$-invariant and thus
$\Lambda\left((P_{\X}A)|_{\X}\right)$ is a subset of
$\Lambda\left((P_{\X+\Y}A)|_{\X+\Y}\right)$, counting the multiplicities, where
the spectrum of $(P_{\X+\Y}A)|_{\X+\Y}$ consists of $\dim(\X+\Y)$ eigenvalues, 
counting the multiplicities, since both $\X$ and $\Y$, and thus their sum $\X+\Y$,
are all finite dimensional.

The only somewhat nontrivial part of the proof is establishing that
the spectrum of the restriction of $(P_{\X+\Y}A)|_{\X+\Y}$
to $\X$ comprises the $\dim\X$ largest eigenvalues of $(P_{\X+\Y}A)|_{\X+\Y}$ using
the lemma assumption that $\X$ is an $A$-invariant subspace
corresponding to the top part of the spectrum of $A$.
In other words, adding $\Y$ to $\X$ does not add any new eigenvalues above
$i=\dim\X$.
We already know that $\Lambda\left((P_{\X}A)|_{\X}\right)$,
on the one hand, makes up the top $\dim\X$ points of the spectrum, which are
eigenvalues, counting the multiplicities, of $A$ and, on the other hand,
is a subset of $\Lambda\left((P_{\X+\Y}A)|_{\X+\Y}\right)$. We only need
to show that the $i=\dim\X$-th eigenvalue of $(P_{\X}A)|_{\X}$, which is
at the same time the $i$-th top point of the spectrum of $A$,
counting the multiplicity of eigenvalues, bounds above the $i+1$-th
eigenvalue of $(P_{\X+\Y}A)|_{\X+\Y}$. But $\Lambda\left((P_{\X+\Y}A)|_{\X+\Y}\right)$
is a vector of Ritz values of $A$ on the trial subspace $\X+\Y$, so this
follows directly from the inf-sup principle
for arbitrary Hermitian (not necessarily compact) operators, see, e.g.,\
\cite[Chapter II, Section 7]{gk} and \cite[Theorem XIII.1]{reedSim}.
\end{proof}

We note that the assumptions of Lemma \ref{l:csle} are not of course applicable 
to all bounded Hermitian operators. E.g.,\  Lemma \ref{l:csle} cannot be applied to 
an orthogonal projector with an infinite dimensional range. It rather covers 
the class of operators with the top part of the spectrum being discrete---a modest, 
but practically important, extension of the class of compact operators; see again
\cite[Chapter II, Section 7]{gk} and \cite[Theorem XIII.1, p. 76]{reedSim}. 
We finally note that the assumption of boundedness (below) of $A$ is not essential 
and can be easily replaced with the assumption that the subspace $\X+\Y$ is in the 
domain of the definition of the corresponding quadratic form.  

The arguments above allow us to substitute the original
infinite dimensional $\H$ and $A$
with finite dimensional $\X+\Y$ and $(P_{\X+\Y}A)|_{\X+\Y}$ in Theorem~\ref{cor1g}.
\begin{theorem} \label{th:inf}
The infinite dimensional, $\dim\H=\infty$, versions of 
Theorem~\ref{cor1g} and its corollary \eqref{eq:SFgn} hold 
under the assumptions of Lemma \ref{l:csle}.
\end{theorem} 

\section{Application to the FEM} \label{s:FEM}
In the FEM context, see, e.g.,\ \cite{MR89k:65132,bo91,MR2206452}, 
let us consider a specific example, the clamped membrane vibration problem---a well 
known eigenvalue problem for the negative Laplacian $-\Delta$ operator 
in two dimensions. 
Let the membrane be a non-convex polygon $\Omega$ with a single reentrant corner
$w\in(\pi,2\pi)$. 
We will use the standard Sobolev spaces $\dot{H}^1(\Omega)$ of functions 
satisfying the homogeneous Dirichlet conditions on the boundary of  $\Omega$ and 
 ${H}^{1+\alpha}(\Omega)$ with $\alpha>0.$

We set $\H=\dot{H}^1(\Omega)$ and define our operator $A$ as, 
informally speaking, the inverse to the negative Laplacian; see, e.g.,\ \cite{MR2206452}, 
so that $A>0$ is compact in $\H$. 
Let us highlight that in this context we use the 
$\H=\dot{H}^1(\Omega)$ scalar product in the definition of the angles to bound the 
largest eigenvalues of $A$, which are the reciprocals 
of the smallest eigenvalues of the negative Laplacian. 
We are looking for an approximation of 
the invariant space $\X\subset\dot{H}^1(\Omega)$ of $A$,
corresponding to the main membrane vibration modes,  
within a trial subspace $\Y\subset\dot{H}^1(\Omega)$ by the Rayleigh-Ritz method.
Using the simplest FEM setup,
the domain $\Omega$ is triangulated according to traditional assumptions, 
and $\Y$ consists of all piecewise linear (on each triangle) continuous functions satisfying 
the homogeneous Dirichlet conditions on the boundary $\partial\Omega$. 
The largest  linear size of the  largest triangle is denoted by $h$.
It holds that $0<{\lmin}_{(\X+\Y)}\to0$ as $h\to0$, so 
we replace ${\lmin}_{(\X+\Y)}$ with its lower bound $0$
in \eqref{eq:SFgn} and Theorem \ref{cor1g}. 

The angles on the right-hand sides in our eigenvalue approximation 
error bounds characterize the approximability
of the target invariant subspace $\X$ by finite element functions from $\Y$, which is 
typically measured by $Ch^\alpha$, where $C$ is a generic constant,
$h$ approaches zero, and the exponent $\alpha$ describes the approximation order. 
The approximability is determined by the type of the FEM, smoothness of 
functions in $\X$, and the choice of the space~$\H$. For our example, 
the approximability bound for a function $v\in{H^{1+\alpha}}(\Omega)$ with some 
$\alpha\in(0,1]$ is 
$\sin\Theta(v,\Y)\leq Ch^\alpha \|v\|_{H^{1+\alpha}}/\|v\|_{H^{1}}$. 
The actual lower bound for $\alpha$, which is $\pi/\omega-\epsilon$, is determined by the angle $\omega$
of the reentrant corner of the polygon  $\Omega$,
which may lead to a corner singularity in eigenfunctions.
The upper bound, $1$, comes from the use of the piecewise linear FEM.  

Let us consider a particular case, where $\dim\X=2$, denoting
the largest eigenvalues by
$\Lambda\left((P_{\X}A)|_{\X}\right)=[\lambda_1,\lambda_2]$ 
and the corresponding $\H$-normalized eigenfunctions by $v_1$ and $v_2$ in $\X$. 
Typically, both eigenfunctions $v_1$ and $v_2$ would have similar 
corner singularities in the reentrant corner, so both $v_1$ and $v_2\in{H^{1+\alpha}}$ 
with $\alpha=\pi/\omega-\epsilon$, but one of their linear combinations, e.g.,\ 
(for illustrative purposes) $v_1-v_2$,  
might have the full $H^2$ regularity, i.e.,\  $v_1-v_2\in{H^{2}}(\Omega)$, and so 
by the approximability result we have $\sin\Theta(v_1-v_2,\Y)\leq Ch$. 
Thus, $\sin\Theta(\X,\Y)\leq C[h^\alpha,h]$; here and below 
we neglect terms that are a smaller order of magnitude in $h$ compared to the terms kept. 
To clarify the example, let us assume that simply $\sin\Theta(\X,\Y)=[h^\alpha,h]$.

This assumption may not be practical for our specific membrane problem. 
However, examples are given in \cite{bo91}, 
where eigenfunctions have different regularities,
while corresponding to the same (multiple) eigenvalue. 
A perturbation argument shows that 
our assumption on regularity of linear combinations of eigenfunctions 
is realistic.

Using the notation
$\Lambda_{\dim{\X}}\left((P_{\Y}A)|_{\Y}\right)=[\lambda_1^h,\lambda_2^h]$
for the relevant FEM Ritz values, we obtain from \eqref{eq:SFgn},  
as in  \cite{MR2206452}, that 
\begin{equation}\label{e:st}
0\leq\lambda_1-\lambda_1^h\leq\lambda_1 h^{2\alpha} \text{ and } 
0\leq\lambda_2-\lambda_2^h\leq\lambda_2 h^{2\alpha},
\end{equation}
while \eqref{eq:1} implies the bound for the error in the trace,
\begin{equation}\label{e:st1}
0\leq\lambda_1+\lambda_2-\lambda_1^h-\lambda_2^h\leq
\lambda_1 h^{2\alpha}+\lambda_2 h^2\approx\lambda_1 h^{2\alpha},
\end{equation}
and \eqref{eq:2} gives the bound for the error in the product,
\begin{equation}\label{e:st2}
0\leq\frac{\lambda_1\lambda_2}{\lambda_1^h\lambda_2^h}-1\lessapprox
(1+h^{2\alpha})(1+h^2)-1\approx h^{2\alpha}.
\end{equation}
The standard bound \eqref{e:st} implies 
\eqref{e:st1} and  \eqref{e:st2} only with an extra factor $2$
in the right-hand side. We conclude that \eqref{eq:SFgn}
cannot take advantage of the
better approximability of the function $v_1-v_2$ in this example,
while our new majorization bounds \eqref{eq:1} and \eqref{eq:2} can, 
and lead to an improvement in the constant with the factor $\dim\X=2$ for the 
trace and product error bounds. 

\section*{Conclusions}

Majorization is a powerful tool that gives elegant and general error bounds
for eigenvalues approximated by the Rayleigh-Ritz method. 
We discover several new results of this kind, including 
multiplicative bounds for relative errors. 
We apply majorization, apparently for first time, in the context of FEM error bounds.
Our initial results are
promising and expected to lead to further development of the majorization
technique for the theory of eigenvalue computations.

\section{Appendix}\label{section_definitions}

Facts on majorization and angles, and most proofs are given here.
\subsection{Weak Majorization}\label{sect.majorization}
For a real vector
$a = \left[ a_{1},\cdots,a_{n} \right]$
let $a^\da$ be %the vector 
obtained by rearranging the entries of $a$ in
an algebraically decreasing order,
$a_{1}^\da\geq\cdots\geq a_{n}^\da.$
We denote $|a|=[|a_1|,\cdots,|a_n|]$ and $a^+=\max\{a,0\}$. 
We say that the vector $b$ weakly majorizes the vector $a$ and we use the notation
$[a_1,\cdots,a_n] \prec_w [b_1,\cdots,b_n]$ or $a\prec_wb$
if $\sum_{i=1}^k a_i^\da\le\sum_{i=1}^k b_i^\da,\quad \ k=1,\ldots,n$.
If in addition the sums above for $k=n$ are equal,
$b$ (strongly) majorizes $a$, which is denoted by  $a\prec b$.
Nonnegative vectors of different sizes may be compared
by appending or removing zeros to match the sizes.

The additive majorization statement $x-y\prec_w z$ for $n$-vectors 
$x=x^\da,\,y=y^\da,$ and $z=z^\da$ is 
equivalent to 
\[\sum_{j=1}^k x_{i_j}\leq\sum_{j=1}^k z_j + y_{i_j},\, 
\forall k:1\leq k\leq n,\, \forall i_j:1\leq i_1<\cdots<i_k \leq n\]
with $x-y\prec z$ if $k=n$ gives the equality. We write  
$\log x-\log y\prec_w\log z$ if 
\[\prod_{j=1}^k x_{i_j}\leq\prod_{j=1}^k z_jy_{i_j},\,
\forall k:1\leq k\leq n,\, \forall i_j:1\leq i_1<\cdots<i_k \leq n,\]
and $\log x-\log y\prec\log z$ if in addition the case $k=n$ gives the equality, 
for nonnegative vectors $x=x^\da,\,y=y^\da,$ and $z=z^\da$. For strictly positive vectors this 
follows directly from the definition of (weak) majorization. 

We need several simple general facts on weak majorization:
If nonnegative vectors $a,b,$ and $c$ are decreasing and of the same size,
then $a\prec_w b$ implies $ac\prec_w bc$, but the converse is not true in general.
If $a\prec_wb\leq c$ then $a\prec_wc$. Concatenation holds, i.e.,\ 
 $a\prec c$ and $b\prec d$ imply $[a,b]\prec[c,d]$; \cite[Corollary II.1.4, p. 31]{bhatia_book}.
If $a\prec_w b$ and $c\prec_w d$ then $a+c\prec_w b+d$
for real vectors, if the bounds $b$ and $d$ are ordered in the same way; \cite[Prop. 4.A.1.b]{mo}.
For a convex increasing function $g(t)$ (e.g.,\ $g(t)=e^t$)
$a\prec_wb$ implies $g(a)\prec_wg(b)$;
\cite[Prop. 4.B.2., p. 109]{mo}. Trivially, $a\prec_wa^+$.  

Let $S(A)$ denote the vector of all singular values of the matrix $A$
in decreasing order; and
for $A$ with real eigenvalues
let $\Lambda(A)$ denote the vector of all eigenvalues of $A$
in decreasing order.
The following theorems are mostly known; see, e.g.,\ \cite{bhatia_book,mo}.
\begin{theorem}[Lidsksi\v{i}]\label{thm:Lid}
$\Lambda(A)-\Lambda(B)\prec \Lambda(A-B)$ for Hermitian $A$ and $B$.
\end{theorem}
\begin{theorem} \label{thm:svp}
$\log S(AB)-\log S(B)\prec\log S(A)$ for general $A$ and $B$, where 
we append zeros to the vectors of singular values if necessary to match the sizes.
\end{theorem}
\begin{proof}
For square matrices (or operators within the same space) 
this is the classical Gelfand-Na\v{i}mark theorem 
\cite[Theorem III.4.5]{bhatia_book}. 
Non-square matrices are extended with zero blocks to obtain square matrices. 
The extension with zero blocks only appends zero singular values and does not change the ranks.
\end{proof}
\begin{theorem}\label{sv_product_majorization}
$S(AB)\prec_w S(A)S(B)$ for general $A$ and $B$, where 
we append zeros to the vectors of singular values if necessary to match the sizes.
\end{theorem}
\begin{proof}
We add $\log S(B)$ to both sides of the statement of Theorem \ref{thm:svp}
and take the exponential function. 
\end{proof}

Our next theorem generalizes Theorem \ref{thm:svp} and 
improves \cite[Corollary 2.4]{li99lidskiimirskywielandt}.
\begin{theorem} \label{thm:svpp}
$\log S(ABC)-\log S(B)\prec\log (S(A)S(C))$ for general $A$, $B$, and $C$,
where we append zeros to singular values if necessary to match the sizes.
\end{theorem}
\begin{proof}
Theorem \ref{thm:svp} can also be formulated as 
$\log S(AB)-\log S(B)\prec\log S(A)$ as singular values of $AB$ and $BA$ 
are the same up to zeros, so Theorem \ref{thm:svp} gives both  
$\log S(ABC)-\log S(BC)\prec\log S(A)$ and $\log S(BC)-\log S(B)\prec\log S(C)$.  
As the right-hand sides in these majorization statements are ordered in the same manner, 
we can add the statements, obtaining the claim of the theorem.
%which in the product form reads as 
%$\prod_{j=1}^k s_{i_j}(ABC) \leq \prod_{j=1}^k s_{i_j}(B)s_j(A)s_j(C)$ for any 
%$k: 1 \leq k \leq n$ and any set of indices $i_j: 1\leq i_1<\cdots<i_k \leq n,$ 
%where $n$ is the largest matrix size and we get the equality if $k=n.$
\end{proof}

We also need the following generalized pinching inequality which may be new.
\begin{theorem}\label{th:pin}
For  matrices  $A_1$, $A_2$, $B$, $C_1$, and $C_2$, 
such that all the products
$A_i^HBC_j$ exist for $i,j=1,2$, we have, possibly up to zeros, 
\begin{equation}\label{eq:genpinching.weak}
[S(A_1^HBC_1),S(A_2^HBC_2)]\prec_w
S\left(\sqrt{A_1A_1^H+A_2A_2^H}B\sqrt{C_1C_1^H+C_2C_2^H}\right),
\end{equation}
and in the case that $A_i=C_i$ and $B=B^H$, we have, possibly up to zeros, 
\begin{equation}\label{eq:genpinching.strong}
[\Lambda(A_1^HBA_1),\Lambda(A_2^HBA_2)]\prec_w
\left(\Lambda\left(\sqrt{A_1A_1^H+A_2A_2^H}B\sqrt{A_1A_1^H+A_2A_2^H}\right)\right)^+.
\end{equation}
\label{thm:genpinching}
\end{theorem}
\begin{proof}
We denote $A=[A_1\,A_2]$ and $C=[C_1\,C_2]$ and form the $2$-by-$2$ block matrix
\[
D=A^HBC=
\left[\begin{array}{cc}
A_1^HBC_1 & A_1^HBC_2 \\
A_2^HBC_1 & A_2^HBC_2 \end{array} \right].
\]
By the standard pinching  inequality, e.g.,\ \cite[Problem II.5.4]{bhatia_book},
the combined singular values of the diagonal blocks of the matrix $D$ are weakly
majorized by the singular values of $D$.
Using the fact that eigenvalues of matrix products
do not depend on the order of the multipliers
shows that the singular values $S(D)$ up to zeros are the same as
$S\left(\sqrt{AA^H}B\sqrt{CC^H}\right)=
S\left(\sqrt{A_1A_1^H+A_2A_2^H}B\sqrt{C_1C_1^H+C_2C_2^H}\right)$,
giving \eqref{eq:genpinching.weak}.

In the case that $A_i=C_i$ and $B=B^H$, the eigenvalues of the
diagonal blocks of $D$ (which are now square) are strongly
majorized by the eigenvalues of $D$. For the latter, we have, up to zeros, 
$\Lambda(D)=\Lambda(A^HBA)=\Lambda(AA^HB)=\Lambda(\sqrt{AA^H}B\sqrt{AA^H})$.
Appending or removing zeros preserve 
majorization for nonnegative vectors, so we replace
$\Lambda(D)\prec_w\left(\Lambda(D)\right)^+$ in the formulas above,
which proves \eqref{eq:genpinching.strong}.
\end{proof}

If $\size[A_2 A_2]=\size B$, there are no zeros appearing, so 
\eqref{eq:genpinching.strong} holds as a strong majorization and without the $+$ operation, 
as in the standard pinching inequality. 

\subsection{Principal Angles Between Subspaces}\label{sect.a}
We need the following:
\begin{theorem}\label{th:a1}\cite[Theorem 3.4]{ka02}
Let $\dim\X=\dim\Y$. Then we have the equalities
$\Lambda\left(P_{\X}P_{\YP}P_{\X}\right)=
S^2\left(P_{\Y}P_{\X^\perp}\right)=S^2\left(P_{\X^\perp}P_{\Y}\right)=
\left[\sin^2\Theta(\X,\Y),0,\ldots,0\right].$
\end{theorem}
\begin{theorem}\label{thm:lambdasin}\cite[Theorem 2.16]{kja07}
If $\dim\X=\dim\Y=p$, then 
the first, i.e.,\ largest, $p$ components of the vector $\Lambda(P_{\X}-P_{\Y})$
are given by the vector $\sin\Theta(\X,\Y).$
\end{theorem}

\subsection{Proofs}\label{s:p}
In this section we provide the main and relatively long proofs.
\subsubsection{Proof of Theorem \ref{thm_proximity_paper}}%
\label{sec:proximity_paper}
We start with two important simplifications.%
\footnote{Ilya Lashuk's proof (private communication, unpublished)
of a particular case, where $A$ is an orthogonal projector,
of bound \eqref{eq.hyp} has stimulated our present proof.}
First, by \cite[Remark 4.1]{ka05} we use the
subspace $\X+\Y$ and the operator $(P_{\X+\Y}A)|_{\X+\Y}$
as substitutions for the original space $\H$ and the original operator $A$,
keeping the same notation, without loss of generality.
Second, the differences of Ritz values do not change with a shift of $A$, i.e.,
for a real $\alpha$ and $A_s=A-\alpha I$, we have, e.g.,\
$\Lambda\left(\left(P_{\X}A\right)|_{\X}\right)-\Lambda\left(\left(P_{\Y}A\right)|_{\Y}\right)=
\Lambda\left(\left(P_{\X}{A_s}\right)|_{\X}\right)-
\Lambda\left(\left(P_{\Y}{A_s}\right)|_{\Y}\right)$.
So all our statements are invariant with respect to a real shift of $A$,
which we can freely choose.
Also, our bounds are invariant with respect to a real scaling of $A$.
Thus for any real $\alpha$ and $\beta\neq0$ we can replace $A$ with
$\beta(A-\alpha I)$ without loss of generality. 
We take  $\alpha=\lmin$ and $\beta=1/(\lmax-\lmin)$
so in the rest of the proof
we assume that $A$ is already shifted and scaled such that
$\lmin=0$ and $\lmax=1$, which guarantees well defined
square roots $\sqrt{A}$ and $\sqrt{I-A}$.

We now prove that
\begin{equation}
|\Lambda((P_{\X}A)|_{\X})-\Lambda((P_{\Y}A)|_{\Y})|\prec_w
\sin \Theta (\X,\Y), \label{eq:mainold}
\end{equation}
and if in addition $\X$ is $A$-invariant then
\begin{equation}
|\Lambda((P_{\X}A)|_{\X})-\Lambda((P_{\Y}A)|_{\Y})|\prec_w
\sin^2 \Theta (\X,\Y). \label{eq:mainnew2}
\end{equation}

Concatenating positive and negative values together, we obtain
\begin{multline*}
\left[\left|\Lambda((P_{\X}A)|_{\X})-\Lambda((P_{\Y}A)|_{\Y})\right|,\,
-\left|\Lambda((P_{\X}A)|_{\X})-\Lambda((P_{\Y}A)|_{\Y})\right|\right]^\da\\
=\left[\Lambda((P_{\X}A)|_{\X})-\Lambda((P_{\Y}A)|_{\Y}),\,
-\left(\Lambda((P_{\X}A)|_{\X})-\Lambda((P_{\Y}A)|_{\Y})\right)\right]^\da\\
=\left[\Lambda((P_{\X}A)|_{\X})-\Lambda((P_{\Y}A)|_{\Y}),\,
\Lambda((P_{\X}(I-A))|_{\X})-\Lambda((P_{\Y}(I-A))|_{\Y})\right]^\da.
\end{multline*}
It is more convenient for us to work in the whole space, so 
above we replace
\begin{multline*}
\left[\Lambda((P_{\X}A)|_{\X})-\Lambda((P_{\Y}A)|_{\Y}),\,
\Lambda((P_{\X}(I-A))|_{\X})-\Lambda((P_{\Y}(I-A))|_{\Y}),0,\ldots,0\right]^\da\\
=\left[\Lambda(P_{\X}AP_{\X})-\Lambda(P_{\Y}AP_{\Y}),\,
\Lambda(P_{\X}(I-A)P_{\X})-\Lambda(P_{\Y}(I-A)P_{\Y})\right]^\da,
\end{multline*}
using
$\left[\Lambda((P_{\X}A)|_{\X}),0,\ldots,0\right]=\Lambda(P_{\X}AP_{\X})$
and similar formulas involving $\Y$ instead of $\X$ and $I-A$ instead of $A$, 
which all hold since $A\geq0$ and $I-A\geq0$ in this proof, 
so the added zeros are correctly placed.

At this point the proof splits for  \eqref{eq:mainold} and \eqref{eq:mainnew2}. 
To prove \eqref{eq:mainold}, in the first sub-vector 
$\Lambda(P_{\X}AP_{\X})-\Lambda(P_{\Y}AP_{\Y})$
we swap  multipliers without changing the eigenvalues, e.g.,\ 
$\Lambda(P_{\X}\sqrt{A}\sqrt{A}P_{\X})=\Lambda\left(\sqrt{A}P_{\X}\sqrt{A}\right)$,
and use Theorem~\ref{thm:Lid},
\begin{eqnarray*}
\Lambda(P_{\X}AP_{\X})-\Lambda(P_{\Y}AP_{\Y})
&=&\Lambda\left(\sqrt{A}P_{\X}\sqrt{A}\right)-
\Lambda\left(\sqrt{A}P_{\Y}\sqrt{A}\right)\\
&\prec&\Lambda\left(\sqrt{A}(P_{\X}-P_{\Y})\sqrt{A}\right);
\end{eqnarray*}
and similarly for the second  sub-vector, 
\[
\Lambda(P_{\X}(I-A)P_{\X})-\Lambda(P_{\Y}(I-A)P_{\Y})
\prec\Lambda\left(\sqrt{I-A}(P_{\X}-P_{\Y})\sqrt{I-A}\right).
\]
%\begin{eqnarray*}
%\Lambda(P_{\X}(I-A)P_{\X})-\Lambda(P_{\Y}(I-A)P_{\Y})
%&=&\Lambda\left(\sqrt{I-A}P_{\X}\sqrt{I-A}\right)-
%\Lambda\left(\sqrt{I-A}P_{\Y}\sqrt{I-A}\right)\\
%&\prec&\Lambda\left(\sqrt{I-A}(P_{\X}-P_{\Y})\sqrt{I-A}\right),
%\end{eqnarray*}
We concatenate, as $a\prec c$ and $b\prec d$ imply $[a,b|\prec[c,d]$,
and we obtain
\[
\begin{split}
[\Lambda(P_{\X}AP_{\X})-\Lambda(P_{\Y}A&P_{\Y}),\,
\Lambda(P_{\X}(I-A)P_{\X})-\Lambda(P_{\Y}(I-A)P_{\Y})]\\
%&=\left[\Lambda\left(\sqrt{A}P_{\X}\sqrt{A}\right)-
%\Lambda\left(\sqrt{A}P_{\Y}\sqrt{A}\right),\,
%\Lambda\left(\sqrt{I-A}P_{\X}\sqrt{I-A}\right)-
%\Lambda\left(\sqrt{I-A}P_{\Y}\sqrt{I-A}\right)
%\right]\\
&\prec\left[\Lambda\left(\sqrt{A}(P_{\X}-P_{\Y})\sqrt{A}\right),\,
\Lambda\left(\sqrt{I-A}(P_{\X}-P_{\Y})\sqrt{I-A}\right)\right]\\
&\prec_w\left[\left(\Lambda(P_{\X}-P_{\Y})\right)^+,0,\ldots,0\right].
\text{~~~~~~~~~~~~~~~(by Theorem \ref{thm:genpinching})}
\end{split}
\]
Picking up the $\dim\X=\dim\Y$ largest (nonnegative) elements on
both sides of this strong majorization statement proves the weak majorization
claim \eqref{eq:mainold}, since by Theorem~\ref{thm:lambdasin}
the $\dim\X=\dim\Y$ largest elements of $\Lambda(P_{\X}-P_{\Y})$ 
are equal to $\sin\Theta (\X,\Y)$.

To prove \eqref{eq:mainnew2} having $A=P_{\X}AP_{\X}+P_{\X^\perp}AP_{\X^\perp},$
we notice that
\begin{equation}\label{eq:less}
\Lambda(P_{\Y}AP_{\Y})
=\Lambda(P_{\Y}P_{\X}AP_{\X}P_{\Y}+P_{\Y}P_{\XP}AP_{\XP}P_{\Y})
\geq \Lambda(P_{\Y}P_{\X}AP_{\X}P_{\Y}),
\end{equation}
as $P_{\Y}P_{\XP}AP_{\XP}P_{\Y}\geq0$, and
similarly $\Lambda(P_{\Y}(I-A)P_{\Y})\geq\Lambda(P_{\Y}P_{\X}(I-A)P_{\X}P_{\Y})$, so
\[
\begin{split}
[\Lambda&(P_{\X}AP_{\X})-\Lambda(P_{\Y}AP_{\Y}),\,
\Lambda(P_{\X}(I-A)P_{\X})-\Lambda(P_{\Y}(I-A)P_{\Y})]\\
&\leq
[\Lambda(P_{\X}AP_{\X})-\Lambda(P_{\Y}P_{\X}AP_{\X}P_{\Y}),\,
\Lambda(P_{\X}(I-A)P_{\X})-\Lambda(P_{\Y}P_{\X}(I-A)P_{\X}P_{\Y})]\\
&\prec
\left[\Lambda\left(\sqrt{A}P_{\X}P_{\YP}P_{\X}\sqrt{A}\right),\,
\Lambda\left(\sqrt{I-A}P_{\X}P_{\YP}P_{\X}\sqrt{I-A}\right)\right]\\
&\prec_w\left[\Lambda(P_{\X}P_{\YP}P_{\X}),0,\ldots,0\right]
\text{~~~~~~~~~~~(by Theorem \ref{thm:genpinching}, %
since $\Lambda(P_{\X}P_{\YP}P_{\X})\geq0$)}\\
&= \left[\sin^2\Theta (\X,\Y),0,\ldots,0\right].
\text{~~~~~~~~~~~~~~~~~~~~~~~~~~~~~~~~~~~~~~~~~~~~~~(by Theorem \ref{th:a1})}
\end{split}
\]
Here, in the second line we again use that
the eigenvalues of the matrix product do not depend on the order
of the matrix multipliers, so we transform, e.g.,\ in the first vector,
$\Lambda(P_{\X}AP_{\X})=\Lambda\left(\sqrt{A}P_{\X}\sqrt{A}\right)$
and $\Lambda(P_{\Y}P_{\X}AP_{\X}P_{\Y})=
\Lambda\left(\sqrt{A}P_{\X}P_{\Y}P_{\X}\sqrt{A}\right)$.
In the next line we independently apply Theorem  \ref{thm:Lid} to each of the two sub-vectors.

\subsubsection{Proof of Theorem \ref{thm_invariant2}}\label{s:mp}
As in the previous proof, we start with two simplifications.
The first one is the same: by \cite[Remark 4.1]{ka05} we use the
subspace $\X+\Y$ and the operator $(P_{\X+\Y}A)|_{\X+\Y}$
as substitutions for the original space $\H$ and the original operator $A$
keeping the same notation, without loss of generality.
Second, we choose
$\alpha=\min\left\{\Lambda\left(\left(P_{\X}A\right)|_{\X}\right)\right\}$ and assume that
the shift is already applied to $A$, i.e.,\ without loss of generality
we assume that both $P_{\X}AP_{\X}$ and  $P_{\XP}(-A)P_{\XP}$ are nonnegative definite 
and so they have well-defined square roots $\sqrt{P_{\X}AP_{\X}}$ 
and $\sqrt{P_{\XP}(-A)P_{\XP}}$, correspondingly. 

For an $A$-invariant subspace $\X$, we split
$A= P_{\X}AP_{\X} + P_{\X^\perp}AP_{\X^\perp},$
and adding and subtracting $\Lambda((P_{\Y}P_{\X}AP_{\X})|_{\Y})$ we derive
\begin{eqnarray*}
0&\leq&\Lambda((P_{\X}A)|_{\X})-\Lambda((P_{\Y}A)|_{\Y})\\
&=&\Lambda((P_{\X}A)|_{\X})
-\Lambda((P_{\Y}P_{\X}AP_{\X})|_{\Y}+(P_{\Y}P_{\XP}A P_{\XP})|_{\Y})\\
&=&\Lambda((P_{\X}A)|_{\X})-\Lambda((P_{\Y}P_{\X}AP_{\X})|_{\Y})\\
& &+\Lambda((P_{\Y}P_{\X}AP_{\X})|_{\Y})
-\Lambda((P_{\Y}P_{\X}AP_{\X})|_{\Y}+(P_{\Y}P_{\XP}AP_{\XP})|_{\Y}).
\end{eqnarray*}

Now, we bound separately the two terms in the sum in the last
two lines. We remind the reader that
$a\prec b$ and $c\prec d$ imply $a+c\prec b^\da+d^\da$
for real vectors,
and this holds similarly for weak majorization.

It is convenient to extend the operators' restrictions by zero to the whole space 
and use a convention that operations
and comparisons of nonnegative decreasing vectors with
different numbers of components is done by appending
zeros at the end of the vectors to match the vectors' sizes, e.g.,\
$\Lambda(P_{\X}AP_{\X})=[\Lambda\left(\left(P_{\X}A\right)|_{\X}\right),0,\ldots,0]\geq0.$
Since $\dim{\X}=\dim{\Y}$ the number of zeros to add for $\X$ and $\Y$
is the same. However, a seemingly trivial claim
$\Lambda(P_{\Y}AP_{\Y})=[\Lambda((P_{\Y}A)|_{\Y}),0,\ldots,0]$
is, in fact, \emph{wrong}, since 
the components of $\Lambda((P_{\Y}A)|_{\Y})$ may not be all nonnegative,
so the added zeros are misplaced compared to $\Lambda(P_{\Y}AP_{\Y}),$ which is 
decreasing by definition. 

We start with the first term in the sum on the right-hand side.
Since both $P_{\X}AP_{\X}\geq0$ and $P_{\Y}P_{\X}AP_{\X}P_{\Y}\geq0$,
we concatenate with zeros correctly and obtain
\begin{eqnarray*}
[\Lambda((P_{\X}A)|_{\X}) -\Lambda((P_{\Y}P_{\X}AP_{\X})|_{\Y}),0,\ldots,0]
&=&\Lambda(P_{\X}AP_{\X}) -\Lambda(P_{\Y}P_{\X}AP_{\X}P_{\Y})\\
&\prec& \Lambda\left(\sqrt{P_{\X}AP_{\X}}P_{\X}P_{\YP}P_{\X}\sqrt{P_{\X}AP_{\X}}\right)\\
&\prec_w& S(P_{\X}AP_{\X})\sin^2\Theta(\X,\Y),
\end{eqnarray*}
applying Theorems \ref{thm:Lid}, \ref{sv_product_majorization}, and \ref{th:a1}.

Considering the second term and again using Theorem \ref{thm:Lid} we get
\[
\begin{split}
\Lambda((P_{\Y}P_{\X}AP_{\X})|_{\Y})-
\Lambda((P_{\Y}P_{\X}AP_{\X})|_{\Y}+(P_{\Y}P_{\XP}AP_{\XP})|_{\Y})\\
\prec\Lambda((P_{\Y}P_{\XP}(-A)P_{\XP})|_{\Y}).
\end{split}
\]
By our assumption on the shift, we have $P_{\XP}(-A)P_{\XP}\geq0$ so
\begin{eqnarray*}
0&\leq&[\Lambda((P_{\Y}P_{\XP}(-A)P_{\XP})|_\Y),0,\ldots,0]\\
&=&\Lambda(P_{\Y}(P_{\XP}(-A) P_{\XP})P_\Y)\\
&=&\Lambda\left(\sqrt{P_{\XP}(-A)P_{\XP}} P_{\XP}P_{\Y}
P_{\XP}\sqrt{P_{\XP}(-A)P_{\XP}}\right)\\
&\prec_w&S(P_{\XP}(-A)P_{\XP})\sin^2\Theta(\X,\Y)
\text{~~~(by Theorems \ref{sv_product_majorization} and \ref{th:a1})}\\
&=&S(P_{\XP}AP_{\XP})\sin^2\Theta(\X,\Y).
\end{eqnarray*}

Adding both bounds together gives the statement of the theorem, i.e.,\
\begin{equation*}
0\leq\Lambda\left(\left(P_{\X}A\right)|_{\X}\right)-
\Lambda\left(\left(P_{\Y}A\right)|_{\Y}\right)
\prec_w \left(S\left(P_{\X}AP_{\X}\right)+S\left(P_{\X^\perp}AP_{\X^\perp}\right)\right)\sin^2\Theta(\X,\Y).
\end{equation*}
%which is equivalent to the statement of the theorem, as
Finally, there are  $\dim(\X+\Y)-\dim\Y$ nonzero components in the vector 
$\Theta(\X,\Y)$, since $\dim(\X+\Y)+\dim(\X\cap\Y)=\dim\X+\dim\Y$.
The first  $\dim(\X+\Y)-\dim\Y$ components in the vector 
$S\left(P_{\X}AP_{\X}\right)+S\left(P_{\X^\perp}AP_{\X^\perp}\right)$
are the same as those in  the vector ${\rm Spr}_{(\X+\Y)}$ 
since we have redefined $A$ such that
the sum $\X+\Y$ gives the whole space and shifted $A$ such that
$\min\left\{\Lambda\left(\left(P_{\X}A\right)|_{\X}\right)\right\}=0$ and so we have 
$S\left(P_{\X}AP_{\X}\right)=\Lambda\left(P_{\X}AP_{\X}\right)$ and
$S\left(P_{\X^\perp}AP_{\X^\perp}\right)=-\Lambda^\ua\left(P_{\X^\perp}AP_{\X^\perp}\right).$
For each component of ${\rm Spr}_{(\X+\Y)}$ where
${\lambda_{\dim\X}}_{(\X+\Y)}<{\lambda_{-i}}_{(\X+\Y)}$
the corresponding angle $\theta_i(\X,\Y)$ in the vector
$\Theta(\X,\Y)$ must be zero, so such a component
of ${\rm Spr}_{(\X+\Y)}$ with an index larger than $\dim(\X+\Y)-\dim\Y$ 
can be defined arbitrarily, since it is multiplied by zero.

\subsubsection{Proof of Theorem \ref{thm:mult}}
First we use exactly the same simplifications as in the beginning of the proof of
Theorem \ref{thm_proximity_paper} in subsection \ref{sec:proximity_paper},
so ${\lmin}_{(\X+\Y)}=0$.
We assume that the space $\H$ is already
mapped into a space of vectors, so that we can use a matrix proof here.
Let $X$ and $Y$ be two matrices
whose columns form orthonormal bases for $\X$ and $\Y$ respectively, so we have
$\Lambda((P_{\X}A)|_{\X})=\Lambda(X^HAX)$ and $\Lambda((P_{\Y}A)|_{\Y})=\Lambda(Y^HAY)$.

The theorem's assumptions $\Lambda((P_{\X}A)|_{\X})>{\lmin}_{(\X+\Y)}$
and $\Theta(\X,\Y)<\pi/2$ give 
$\Lambda((P_{\X}A)|_{\X})\geq\Lambda((P_{\Y}A)|_{\Y})>{\lmin}_{(\X+\Y)}$
by  \eqref{eq:SF}. This is equivalent in our simplified situation to
$\Lambda(X^HAX)\geq\Lambda(Y^HAY)>0$, so we can legitimately 
take the log of their ratio below.  
By analogy with \eqref{eq:less}, since $\X$ is $A$-invariant and
$P_{\YP}AP_{\YP}\geq0$ because of the shift of A that made $A\ge 0$, we have
\begin{equation*}%\label{eq:relative1}
\Lambda(Y^HAY)\geq\Lambda(Y^HP_{\X}AP_{\X}Y)=\Lambda((Y^HX)X^HAX(X^HY))=
\Lambda(C^HX^HAXC),
\end{equation*}
where we denote $C=X^HY$. 
%Since $\min\X=\dim\Y$, 
We have %$\Lambda(C^HC)=\Lambda(CC^H)$, but 
$\Lambda^\ua(CC^H)=\cos^2\Theta(\X,\Y)>0$
by definition and the theorem's assumption, so  
both matrices $C$ and $C^H$ are invertible and then
$1\leq\Lambda^\ua(C^{-H}C^{-1})=\cos^{-2}\Theta(\X,\Y)$. 
The key step is using Theorem \ref{thm:svp}, substituting 
$A:=C^{-H}$ and $B:=C^H\sqrt{X^HAX}$ in
\begin{eqnarray}\nonumber
0\leq\log\left(\frac{\Lambda(X^HAX)}{\Lambda(C^HX^HAXC)}\right)
&=&2\log\left(\frac{S\left(\sqrt{X^HAX}\right)}{S\left(C^H\sqrt{X^HAX}\right)}\right)\\\nonumber
&=&2(\log S(AB)-\log S(B))\\\label{e:log}
&\prec&2\log S(A)\\\nonumber 
&=&\log\Lambda(C^{-H}C^{-1})\\\nonumber
&=&\log\left({\cos^{-2}\Theta(\X,\Y)}\right))\\\nonumber
&=&\log\left({1+\tan^{2}\Theta(\X,\Y)}\right).
\end{eqnarray}
%$B=X^HAX$ and $C$ just defined,
%which gives
%\begin{eqnarray*}
%0\leq\log\left(\frac{\Lambda(X^HAX)}{\Lambda(C^HX^HAXC)}\right)
%&=&\log\left|\frac{\Lambda(X^HAX)}{\Lambda(C^HX^HAXC)}\right|\\
%&\prec_w&\log|\Lambda(C^HC)|\\
%&=&\log\left|\cos^2\Theta(\X,\Y)\right|\\
%&=&\log\left({\cos^{-2}\Theta(\X,\Y)}\right).
%\end{eqnarray*}
Replacing here $\Lambda(C^HX^HAXC)$ with 
$\Lambda(Y^HAY)\geq\Lambda(C^HX^HAXC)$, as shown above, gives the 
multiplicative weak majorization bound of Theorem \ref{thm:mult}.

If $x\prec y$ then $\phi(x)\prec_w\phi(y)$ for
any nondecreasing convex real valued function $\phi$, see,
e.g.,\ \cite[Statement 4.B.2]{mo}.
Taking $\phi(t)=e^t$ for \eqref{e:log} gives
\begin{equation*}
1\leq\frac{\Lambda(X^HAX)}{\Lambda(Y^HAY)} \leq
\frac{\Lambda(X^HAX)}{\Lambda(Y^HP_{\X}AP_{\X}Y)}
=\frac{\Lambda(X^HAX)}{\Lambda(C^HX^HAXC)}
\prec_w 1+\tan^{2}\Theta(\X,\Y).
\end{equation*}
Subtracting the vector of ones gives the second bound of Theorem \ref{thm:mult}.

\section*{Acknowledgments}
The authors thank: 
Mark Embree and anonymous referees for their useful and detailed comments; 
Chris Paige and Ivo Panayotov for stimulating collaboration in \cite{akpp06}, where 
bound \eqref{eq.hyp} has been proved for extreme eigenvalues; 
Ilya~Lashuk for sharing his unpublished proof of a particular case of bound \eqref{eq.hyp} 
with $A=A^2$; 
and Peizhen Zhu for proofreading the manuscript. 

\vskip12pt

%\newpage 

%\bibliographystyle{plainnat}
\def\refname{
\end{document}